\documentclass{amsart}
\usepackage{color,graphicx,wrapfig}
\title[A parabolic free boundary problem]{A parabolic free boundary problem with Bernoulli type condition on the free boundary}
\author[J. Andersson]{John Andersson}
\address{Max Planck Institute for Mathematics in the Sciences, Inselstr. 22, D-04103 Leipzig, Germany}
\email{anders@mis.mpg.de}
\author[G.S. Weiss]{Georg S. Weiss}
\address{Graduate School of Mathematical Sciences,
University of Tokyo, 3-8-1 Komaba, Meguro-ku, Tokyo-to, 153-8914 Japan,}
\email{gw@ms.u-tokyo.ac.jp}
\urladdr{http://www.ms.u-tokyo.ac.jp/~gw/}
\thanks{$2000$ {\it Mathematics Subject Classification.\/} Primary
35R35, Secondary 35K55.}
\thanks{J. Andersson has been partially supported by a fellowship
of the Max Planck Society.
G.S. Weiss has been partially supported by the Grant-in-Aid
15740100/18740086 of the Japanese Ministry of Education, Culture, Sports, Science and Technology and partially supported
by a fellowship of the Max Planck Society. Both authors thank the Max Planck
Institute for Mathematics in the Sciences for the hospitality
during their stay in Leipzig.}
\date{}

\theoremstyle{plain}
\newtheorem{theorem}{Theorem}[section]
\newtheorem{lemma}[theorem]{Lemma}
\newtheorem{proposition}[theorem]{Proposition}
\newtheorem{corollary}[theorem]{Corollary} 
\theoremstyle{definition}
\newtheorem{definition}[theorem]{Definition} 
\theoremstyle{example}

\theoremstyle{definition}
\newtheorem{remark}[theorem]{Remark}
\numberwithin{equation}{section}

\def\R{{\bf R}}

\def\N{{\bf N}}
\def\C{{\bf C}}

\def\sr{{\scriptsize\bf R}}
\def\H{{\bf H}}

\def\div{\textrm{\rm div }}

\def\pardist{\textrm{\rm pardist}}

\def\supp{\textrm{\rm\small supp }}

\begin{document}
\begin{abstract}
Consider the parabolic free boundary problem
$$
\Delta u - \partial_t u = 0 \textrm{ in } \{ u>0\}\> , \>
|\nabla u|=1 \textrm{ on } \partial\{ u>0\}\; .
$$
For a realistic class of solutions, containing for example
{\em all} limits of the singular perturbation problem
$$\Delta u_\varepsilon - \partial_t u_\varepsilon = \beta_\varepsilon(u_\varepsilon)\; \textrm{ as } \varepsilon\to 0,$$
we prove that one-sided flatness of the free boundary implies regularity.\\
In particular, we show that the topological
free boundary $\partial\{ u>0\}$ can be decomposed
into an {\em open} regular set (relative to $\partial\{ u>0\}$)
which is locally a surface with H\"older-continuous
space normal, and a closed singular set.\\
Our result extends the main theorem in the
paper by H.W. Alt-L.A. Caffarelli (1981) to more general solutions
as well as the time-dependent
case. Our proof uses methods developed in
H.W. Alt-L.A. Caffarelli (1981), however we replace the core
of that paper, which relies on non-positive mean curvature
at singular points, by an argument based on scaling
discrepancies, which promises to be applicable to
more general free boundary or free discontinuity problems. 
\end{abstract}
\maketitle
\section{Introduction}
\begin{figure}[hb]
\begin{center}
\input{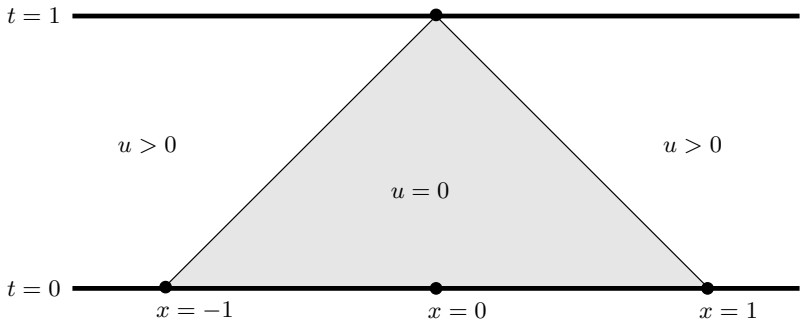}
\end{center}
\caption{Colliding traveling waves}\label{travel}
\end{figure}

\begin{figure}[b]
\begin{center}
\input{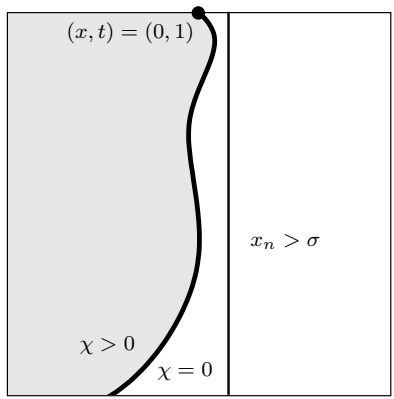}
\end{center}
\caption{One-sided flatness in the case $\rho=1$}\label{onesided}
\end{figure}

\begin{figure}[t]
\begin{center}
\input{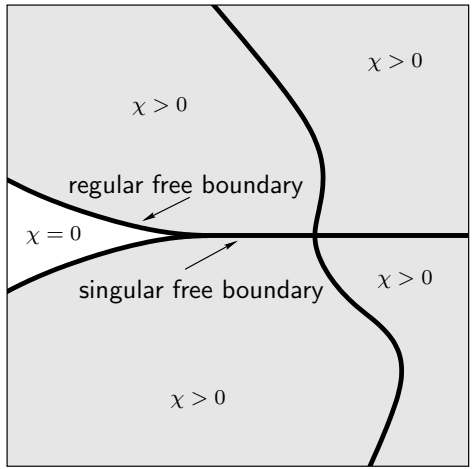}
\end{center}
\caption{Example of the set of regular free boundary points (stationary)}\label{example}
\end{figure}

The parabolic free boundary problem
\begin{equation}\label{bernoulli}\Delta u - \partial_t u = 0 \textrm{ in } \{ u>0\}\> , \>
|\nabla u|=1 \textrm{ on } \partial\{ u>0\}\end{equation}
has originally been derived as singular limit
from a model for the propagation of equidiffusional premixed
flames with high activation energy (\cite{buckmaster});
here $u=\lambda(T_c-T)\> ,$ $T_c$ is the flame
temperature, which is assumed to be constant, $T$ is
the temperature outside the flame and $\lambda$ is a 
normalization factor.
\\
Let us shortly summarize the mathematical results
directly relevant in this context, beginning with
the limit problem (\ref{bernoulli}):
in the brilliant paper \cite{AC}, H.W. Alt and L.A. Caffarelli
proved via minimization of the energy $\int (\vert\nabla u\vert^2
\> + \> \chi_{\{u>0\}})$ --
here $\chi_{\{u>0\}}$ denotes the characteristic function 
of the set ${\{u>0\}}$ --
existence of a stationary solution
of (\ref{bernoulli}) in the sense of distributions.
They also derived regularity of the free boundary
$\partial\{u>0\}$ up to a set of vanishing ${n-1}$-dimensional
Hausdorff measure. 
By \cite{cpde} existence of singular minimizers implies the existence of
singular minimizing cones. 
L.A. Caffarelli-D. Jerison-C. Kenig showed that singular minimizing cones
do not exist in dimension $3$ (\cite{david1}).
Moreover it is known that singular minimizing cones exist
for $n\ge 7$ (\cite{david2}).
{\em Non-minimizing} singular cones appear already
for $n= 3$ (see \cite[example 2.7]{AC}). Moreover it is known,
that solutions of the Dirichlet problem in two space dimensions
are not unique (see
\cite[example 2.6]{AC}).
\\
For the time-dependent (\ref{bernoulli}), both ``trivial non-uniqueness''
(the positive solution of
the heat equation is always another solution of (\ref{bernoulli}))
and ``non-trivial uniqueness'' (see \cite{nonunique}) occur.
Even for flawless initial data, classical solutions of
(\ref{bernoulli}) develop singularities after a finite time
span; consider e.g. the example of two colliding traveling waves
\begin{equation}
\begin{array}{ll}\label{ex1}
u(t,x) = & \chi_{\{x+t>1\}} (\exp(x+t-1)-1)
\\ 
& \quad + \> \chi_{\{-x+t>1\}} (\exp(-x+t-1)-1) \textrm{ for } t \in [0,1)
\end{array}
\end{equation}
(see Figure \ref{travel}).
\\[1cm]
There are several approaches concerning the construction
of a solution of the time-dependent problem, all of which
are based in some form on the convergence of
the solution $u_\varepsilon$ of the reaction-diffusion equation
\begin{equation}\label{perturb}
\Delta u_\varepsilon - \partial_t u_\varepsilon = \beta_\varepsilon(u_\varepsilon)
\end{equation}
to (\ref{bernoulli}) as $\varepsilon\to 0$;
here $\beta_\varepsilon(z)
= {1\over \varepsilon} \beta({z\over \varepsilon})\> , \>
\beta\in C^1_0([0,1])\> , \> \beta > 0$ in $(0,1)$ and
$\int \beta = {1\over 2}\> .$\\
L.A. Caffarelli and J.L. Vazquez 
proved in \cite{cava}
uniform estimates for (\ref{perturb})
and a convergence result: for initial data $u^0$
that are strictly mean concave in the interior of their 
support, a sequence of $\varepsilon$-solutions converges
to a solution of (\ref{bernoulli}) in the sense of
distributions.\\
Let us also mention several results on the corresponding 
two-phase problem, which are relevant as solutions of
the one-phase problem are automatically solutions of the
corresponding two-phase problem. In \cite{leder1} and \cite{leder2},
L.A. Caffarelli, C. Lederman and N.
Wolanski prove convergence to a barrier solution
in the case that the limit function satisfies $\{u=0\}^\circ=\emptyset\> .$
\\
Then, there is the convergence to a solution in the sense of domain
variations \cite{calc} which seems to contain more information
than the barrier solutions in \cite{leder1} and \cite{leder2}.
For more general two-phase problems see \cite{ejde}.
Domain variation solutions play an important rule in this paper and
will be discussed in more detail in Section \ref{notion}.
\\[.3cm]
Here let it suffice to say that domain variation solutions
are pairs $(u,\chi)$ where the order parameter $\chi$
shares many properties with the characteristic function
$\chi_{\{ u>0\}}$ but does not necessarily coincide with
it. By \cite{calc}, {\em all limits} of the singular
perturbation problem (\ref{perturb}) are domain variation solutions,
so all results in the present paper hold for all limits of (\ref{perturb}).\\
Our main result Theorem \ref{theo:main} states
-- leaving out inessential assumptions --
that if
$(0,\rho^2)$ is a point on the topological free boundary
and if the set $\{ \chi >0\}$ is flat enough, i.e.
\begin{displaymath}
\chi(x,t)=0 \textrm{ when } (x,t)\in Q_\rho \textrm{ and }x_n\ge \sigma \rho, 
\end{displaymath}
for some $\sigma\le \sigma_0$ (see Figure \ref{onesided}), then
the free boundary
$Q_{\rho/4} \cap \partial\{ u>0\}$ is
a surface with H\"older-continuous
space normal.\\
As a consequence we obtain that the regular set is open relative
to $\partial\{ u>0\}$ (Corollary \ref{regular}, cf. Figure \ref{example}).
\\
Note that even in the stationary case our result
extends the result in \cite{AC}
as our assumptions {\em do not exclude}
degenerate points or cusps close to the origin
(excluded by the definition of weak solutions \cite[5.1]{AC}),
{\em our result does that}.
\\
In the proof of our result we use ingenious tools developed in
 \cite{AC}:
We prove that flatness on the side of $\{ \chi=0\}$
implies flatness on the side of $\{ \chi>0\}$ which in turn yields
uniform convergence of an inhomogeneously
scaled sequence of free boundaries.
\\
However we replace
the core in the
method of H.W. Alt-L.A. Caffarelli, relying on non-positive
mean curvature of $\partial\{ u>0\}$ at singularities, 
by a method based on {\em scaling discrepancies}
(Proposition \ref{discr}). This original component
gives hope that the method may now be applicable
to more general free boundary or free discontinuity
problems, in particular two-phase free boundary problems.
\section{Notation}
Throughout this article $\R^n$ will be equipped with the Euclidean
inner product $x\cdot y$ and the induced norm $\vert x \vert\> ,\>
B_r(x_0)$ will denote the open $n$-dimensional ball of center
$x_0$, radius $r$ and volume $r^n\> \omega_n\> ,
\> B'_r(0)$ the open $n-1$-dimensional ball of center $0$ and radius
$r$,
and $e_i$ the $i$-th unit vector in $\R^n$.
We define $Q_r(x_0,t_0) := B_r(x_0)\times (t_0-r^2 , t_0+ r^2)$ to be 
the cylinder
of radius $r$ and height $2r^2$,
$Q^-_r(x_0,t_0) := B_r(x_0)\times (t_0-r^2 , t_0)$ its ``negative part''
and
$T_r^-(t_0) := \R^n \times (t_0 - 4r^2 , t_0 - r^2)$ 
the horizontal layer from $t_0 - 4r^2$ to $t_0-r^2$.
Let us also introduce the parabolic distance
$\pardist((t,x),A) := \inf_{(s,y)\in A}
\sqrt{\vert x-y\vert^2 + \vert t-s\vert}$.
Considering 
a function $\phi \in H^{1,2}_{\rm loc}(\R^n;\R^n)$ 
we denote by $\div\phi := \sum_{i=1}^n \partial_i \phi_i$
the space
divergence and
by \[ D\phi := \left(\begin{array}{ccc}
\partial_1 \phi_1 & \dots & \partial_n \phi_1\\
& \dots & \\
\partial_1 \phi_n & \dots & \partial_n \phi_n\end{array}\right)\]
the matrix of the spatial partial derivatives.\\
Given a set $A\subset \R^n\> ,$ we denote its interior by $A^\circ$
and its characteristic function by $\chi_A\> .$
In the text we use the $n$-dimensional Lebesgue-measure
${\mathcal L}^n$ and
the $m$-dimensional Hausdorff measure
${\mathcal H}^m$.
When considering a given set $A\subset \R^n$,
let \[\partial_M A := \{x\in \R^n\> : \>
\limsup_{r\to 0} {{\mathcal L}^n(B_r(x)\cap A)\over
{{\mathcal L}^n(B_r)}} > 0 \textrm{ and } 
 \limsup_{r\to 0} {{\mathcal L}^n(B_r(x) - A)\over
{{\mathcal L}^n(B_r)}} > 0 \}\] be
the measure-theoretic boundary of $A$, let
$\partial^* A := \{x\in \R^n\> : \>
\textrm{ there is } \nu(x) \in \partial B_1(0)
\textrm{ such that }
r^{-n} \int_{B_r(x)} \vert \chi_A-\chi_{\{y:(y-x)\cdot \nu(x)<0\}}
\vert \to 0 \textrm{ as } r\to 0\}$
(by \cite[Corollary 5.6.8]{ziemer} $\partial^* A$ coincides
${\mathcal H}^{n-1}$-a.e. with the reduced boundary of a set
of finite perimeter defined in \cite[Definition 5.5.1]{ziemer}),
and let
$\nu:\partial^* A\to \partial B_1(0)$ denote this measure theoretic
outward normal to $\partial A$.
We shall often use
abbreviations for inverse images like $\{u>0\} := 
\{x\in \Omega\> : \> u(x)>0\}\> , \> \{x_n>0\} := 
\{x \in \R^n \> : \> x_n > 0\}\> , \>
\{s=t\} := \{(s,y)\in \R^{n+1}\> : \> s=t\}$ etc. 
as well as $A(t) := A \cap \{s=t\}$ for a set $A\subset \R^{n+1}$,
and occasionally 
we employ the decomposition $x=(x',x_n)$ of a vector $x\in \R^n$
as well as the corresponding decompositions of the gradient
and the Laplace operator,
$$ \nabla u = (\nabla' u,\partial_n u) \textrm{ and } \Delta u = \Delta' u\> + \> \partial_{nn} u\; .$$
\\
Finally, $\C^{\beta,\mu}:=\H^{\mu,\beta}$
denotes the parabolic
H\"older-space defined in 
\cite{lady}.
\section{Notion of solution and Preliminaries}\label{notion}
In this section we gather some results from \cite{calc}. As degenerate points are
unavoidable in the parabolic problem (see the introduction of \cite{calc}
for examples), an extension of the {\em weak solutions} in \cite{AC}
does not seem to be the right choice. Instead we use the solutions of
\cite[Definition 6.1]{calc}, which are, roughly speaking, solutions in the
sense of domain variations. The advantage is that the class of
solutions defined in \cite[Definition 6.1]{calc} is closed under the blow-up
process. Moreover, {\em all} limits of the singular perturbation problem
discussed in \cite{cava} 
{\em are} domain variation solutions and
satisfy \cite[Definition 6.1]{calc} (see \cite[Section 6]{calc}).
Let us recall the definition of solutions
and the monotonicity formula used therein:
\begin{theorem}[Monotonicity Formula, \textrm{cf. \cite[Theorem 5.2]{calc}}]
\label{mon}
Let $(x_0,t_0)\in \R^n\times (0,\infty)\> , \>
T_r^-(t_0)$ $= \R^n \times (t_0 - 4r^2 , t_0 - r^2)\> ,$
$0<\rho<\sigma<{\sqrt{t_0}\over 2}$ and
\[ G_{(x_0,t_0)}(x,t) =
4\pi (t_0-t)\> {\vert 4\pi (t_0-t)\vert}^{-{n \over 2}-1} 
\; \exp\left(-{\vert x-x_0\vert^2 \over
{4 (t_0-t)}}\right)\; \; .\]
Then
\[ \begin{array}{l}
\Psi_{(x_0,t_0)}(r)\; 
=\; 
r^{-2} \int_{T_r^-(t_0)} \left( 
{\vert \nabla u \vert}^2 \> +\>
\chi
\right)\> G_{(x_0,t_0)}
- \; {1 \over 2} \> r^{-{2}} 
\int_{T_r^-(t_0)}
{1 \over {t_0-t}} \> u^2 \> G_{(x_0,t_0)}
\end{array}\]
satisfies the monotonicity formula 
\[ \Psi_{(x_0,t_0)}(\sigma)\> - \> \Psi_{(x_0,t_0)}(\rho)
\]\[ \ge \;
\int_\rho^\sigma r^{-1-{2}} \>
\int_{T_r^-(t_0)} {1 \over {t_0-t}}\> \Bigg(\nabla u
 \cdot (x-x_0) \> - \>
2 (t_0-t) \partial_t u \> - \> 
u \Bigg)^2 \> G_{(x_0,t_0)}\> dr\; \ge \; 0\; \; .
\]
\end{theorem}
\begin{definition}[\textrm{cf. \cite[Definition 6.1]{calc}}]\label{solution}
We call $(u,\chi)$ a solution in $\Omega_0 := \R^n\times (0,\infty)$
(in which case we set $\tau := 0$) or $\Omega_1 := \R^n\times (-\infty,\infty)$
(in which case we set $\tau := 1$), if:\\
1) $u \in \C^{1,{1\over 2}}_{\rm loc}(\Omega_\tau)\cap
C^2(\Omega_\tau\cap \{u>0\})\cap H^{1,2}_{\rm loc}(\Omega_\tau)$ and
$\chi\in L^1((-\tau R,R);BV(B_R(0)))$ for each $R\in (0,\infty)\> .$
For each $R\in (0,\infty)$ and $\delta\in (0,1)$ there exists
$C_1<\infty$ such that for $Q_r(x_0,t_0)\subset \Omega_\tau\cap Q_R(0)$
\[ \int_{Q_r(x_0,t_0)} \vert \nabla\chi\vert \> \le \> C_1 \> r^{n+1},\]
\[ \int_{Q_r(x_0,t_0)} \vert \partial_t u\vert^2 \> \le \> C_1 \> r^{n}, \textrm{ and}\]
\[ \int_{B_r(x_0)\times (t_0+S_1r^2,t_0+S_2r^2)} \vert\partial_t (\vert\nabla u\vert^2
\> + \> \chi)*\phi_{r\delta} \vert \; \le \; C_1 \sqrt{
S_2-S_1}\> r^n\]
for $0< S_1 < S_2 < \infty$;
here the mollifier $(\phi_\delta)_{\delta\in (0,1)}$ should be non-negative
and satisfy
$\phi_\delta(\cdot) = {1\over {\delta^n}}\phi({\cdot\over {\delta}}),
\phi\in C^{0,1}_0(\R^n)\> ,$ $\int \phi = 1$ and
$\supp \phi \subset B_1(0)\> .$
\\
Moreover, $\chi\in \{0,1\}$ a.e. in $\Omega_\tau$ and 
$\chi_{\{u>0\}} \le \chi$ a.e. in $\Omega_\tau\> .$\\
2) The solution $u$ satisfies the monotonicity formula Theorem \ref{mon}
 (in the case of $\tau=1$ for $(x_0,t_0)\in \R^{n+1}$
and $\sigma\in (0,\infty)$).
\\
\[ \textrm{3) } 0 \; = \; 
\int_{-\infty}^{\infty} \int_{\sr^n} [-2\partial_t u \> \nabla u\cdot \xi
\> + \> (\vert\nabla u\vert^2 \> + \> \chi)\> \div \xi
\> - \> 2 \nabla u D\xi \nabla u]\]
for every $\xi\in C^{0,1}_0(\Omega_\tau;\R^n)\> .$
\\
4) The solution $u$ is non-negative.\\
5) The solution $u$ attains the initial data $u^0\in C^{0,1}_0(\R^n)$
in $L^2_{\rm loc}(\R^n)$ in the case that $\tau=0\> .$
\\
6) For each $\kappa>0$ there is $\delta>0$ such
that $Q_r(x_0,t_0)\subset \Omega_\tau$ and
$\Vert{u(x_0+rx,t_0+r^2 t)\over r}-\theta \vert x_n \vert\Vert_{
C^0(Q_1(0))} < \delta$ imply
$\theta < 1 + \kappa\> .$
\\
7) For $\delta\in (0,1)\> , \>
\psi_\delta \in C^{0,1}_0(\{ \vert y \vert^2 + s^2 < \delta^2
\})\> , \> u_r(y,s) := {u(t_0+r^2 s,x_0+ry)\over r}$ and
$\chi_r(y,s) := \chi(x_0+ry,t_0+r^2 s)$ the following holds:
\[ \textrm{a) }
\; \int_{Q_\rho(x_1,t_1)} \vert (\nabla\chi_r
\cdot x \> + \> 2t\partial_t \chi_r)*\psi_\delta\vert
\]\[ \le \; C(\delta,Z,T,S,\rho)
\left(\Psi_{(x_0,t_0)}(r
\sqrt{{-t_1+\delta+\rho^2\over 2}}) 
\> - \> \Psi_{(x_0,t_0)}(r
\sqrt{{-t_1-\delta-\rho^2\over 2}})\right)\]
for $-S\le t_1 \le -T<0\> , \> \delta+\rho^2\le {T\over 2}\> , \>
\vert x_1 \vert \le Z$ 
and, in the case of $\tau = 0\> ,\> t_0-2r^2(-t_1+\rho^2+\delta)>0\> .$
\[ \textrm{b) }
\; \int_{Q_\rho(t_1,x_1)} \vert (\nabla \chi_r
\cdot \xi) * \psi_\delta\vert
\; \le \; C(\delta) \int_{Q_{\sqrt{\delta}+\rho}(t_1,x_1)} 
\vert \nabla u_{r}\cdot \xi\vert \]
for $\xi \in \partial B_1(0)\> , \> t_1<0$ and, in the case of $\tau = 0 
\> , \> t_0-r^2(-t_1+(\rho+\sqrt{\delta})^2)>0\> .$
\[ \textrm{c) }
\; \int_{t_1}^{t_2} \partial_t((\vert \nabla u_r\vert^2
+\chi_r)
*\phi_\delta)(t,x_0)
\; \le \; \int_{t_1}^{t_2} \int_{\sr} 2\partial_t u_{r}(t,z)
\nabla u_{r}(t,z)\cdot\nabla \phi_\delta(x_0-z)\> dz\]
for $-\infty<t_1<t_2<\infty$
and, in the case of $\tau = 0\> , \> t_0+r^2t_1>0\> .$
\end{definition}
\begin{remark}
As the function $\chi$ is defined only almost everywhere, all
pointwise equalities/inequalities involving $\chi$ should
be understood as equalities/inequalities that hold almost everywhere
with respect to the Lebesgue measure.
\end{remark}
The reader may wonder whether a solution in the sense of distributions
(possibly defined by the identity in \cite[Lemma 11.3]{calc})
would not be good enough for the purposes of this paper.
It turns however out that the information yielded by the order parameter $\chi$ in Definition \ref{solution}
carries information that is essential in what follows.
Incidentally, $\chi$ may be different from $\chi_{\{ u>0\}}$ (see \cite[Remark 4.1]{calc}).
\begin{lemma}\label{prem}
Let $(u,\chi)$ be a solution in the sense of Definition \ref{solution} and suppose that
for some $(x_0,t_0)$ in the set of definition and for
some sequence $r_m\to 0,m\to \infty$
\[ u_{r_m}(y,s) := {u(x_0+{r_m}y,t_0+{r_m}^2 s)\over {r_m}} \to 0
\textrm{ locally in } \{ y_n<0\}\times (-\infty,0) \textrm{ as } m\to\infty\]
and
\[ \chi_{r_m}(y,s) := \chi(x_0+{r_m}y,t_0+{r_m}^2 s)
\to 0 \textrm{ a.e. in }  \{ y_n >0\}\times(-\infty,0)\textrm{ as } m\to\infty\; .\]
Then for some $\delta>0$, $u$ is caloric in
$Q_\delta(x_0,t_0)$ and satisfies
\[ u=0 \textrm{ in } Q^-_\delta(x_0,t_0)\; .\]
\end{lemma}
\proof
The assumptions imply 
by Definition \ref{solution} 1) that
\[ u_{r_m}
\to 0 \textrm{ a.e. in }  \{ x_n >0\}\times(-\infty,0)\textrm{ as } m\to\infty\; .\]
Moreover, they imply by \cite[Proposition 10.1 2)]{calc} that the density
\[ \Psi_{(x_0,t_0)}(0+) \in \{ 0 \} \cup \{ H_n \}\; ,\]
where $H_n$ is the energy of the half-plane solution
defined in \cite[Section 10]{calc}.
In the case 
\[ \Psi_{(x_0,t_0)}(0+)=0\]
we obtain from \cite[Proposition 10.1 2)]{calc} immediately
the statement of the lemma.\\
In the case \[ \Psi_{(x_0,t_0)}(0+)=H_n\]
it follows from  \cite[Proposition 10.1 1)]{calc} that
the limit of $u_{r_m}(y,s)$ as $m\to \infty$ must after rotation be the half-plane
solution $\max(-x_n,0)$, a contradiction to the limit of $u_{r_m}$ being $0$.
\qed
\section{Flatness Classes}
\begin{definition}\label{flatness}
Let $0<\sigma_+, \sigma_- < 1$ and $\tau\ge 0$. We say that 
\begin{displaymath}
u\in F(\sigma_+,\sigma_-,\tau) \quad \textrm{in } Q_{\rho}\quad \textrm{in direction }e_n
\end{displaymath}
if
\begin{enumerate}
\item $(u,\chi)$ is a solution in the sense of Definition \ref{solution} in a domain containing $Q_{\rho}$.
\item 
\begin{displaymath} 
(0,\rho^2) \in \partial \{u>0\},
\end{displaymath}
\begin{displaymath}
u(x,t)=\chi(x,t)=0 \textrm{ when } (x,t)\in Q_\rho \textrm{ and }x_n\ge \sigma_+ \rho, 
\end{displaymath}
\begin{displaymath}
\chi(x,t)=1  \textrm{ and }
u(x,t)\ge -(x_n+\sigma_- \rho) \textrm{ when } (x,t)\in Q_\rho \textrm{ and }x_n\le -\sigma_- \rho\; .
\end{displaymath}
\item 
\begin{displaymath}
|\nabla u|\le 1+\tau \textrm{ in } Q_{\rho}.
\end{displaymath}
\end{enumerate}
When the origin is replaced by $(x_0,t_0)$ and the flatness direction
$e_n$ is replaced by $\nu$ then we define $u$ to belong to the flatness
class $F(\sigma_+,\sigma_-,\tau) $ in $Q_{\rho}(x_0,t_0)$ in direction $\nu$.
\end{definition}
\section{Flatness on the side of $\{ \chi=0\}$ implies
flatness on the side of $\{ \chi>0\}$}
The aim of this and the following sections is to draw information from
properties of an inhomogeneous blow-up limit. One of the central problems when
using blow-up arguments is ``{\em not-strong convergence}'' or
``{\em energy loss}'' in the limit. Here we avoid
those problems by working with {\em uniform convergence} (not
some Sobolev norm). The approach is based on a powerful idea by H.W. Alt-L.A. Caffarelli
who used ``flatness on the side of $\{ u=0\}$ implies
flatness on the side of $\{ u>0\}$'' to prove uniform convergence
to an inhomogeneous blow-up limit (cf \cite[Section 7]{AC}). In this section we extend their
result to a
weaker class of solutions and to
the parabolic case, using results in \cite{calc}.\\
The following Lemma is the parabolic version of \cite[Lemma 4.10]{AC}.
\begin{lemma}\label{lem:ball}
Let $(u,\chi)$ be a solution in the sense of Definition \ref{solution} in a domain containing
the closure of a non-empty open ball $B = \{ (y,s): |(y,s)-(y_0,s_0)| < c\}$
such that
$B\subset \{ \chi=0\}$ and
$B$ touches the set $\{ u>0\}$ at the origin.\\
Then
\begin{displaymath}
\limsup_{\{ u>0\}\ni (x,t)\rightarrow 0}\frac{u(x,t)}{\pardist((x,t),B)} = 1.
\end{displaymath}
\end{lemma}
\proof
Let $Y_k=(y_k,s_k)\in \R^{n+1}$ be a sequence such that
\begin{displaymath}
\ell=\limsup_{\{ u>0\}\ni (x,t)\rightarrow 0}\frac{u(x,t)}{\pardist((x,t),B)}=
\lim_{k\rightarrow \infty}\frac{u(Y_k)}{\pardist(Y_k,B)}.
\end{displaymath}
Set $d_k :=\pardist(Y_k,B)$ and let $(x_k,t_k)=X_k\in \partial B$ be such that 
$\pardist(Y_k,X_k)=d_k$.
\\
We consider the blow-up sequence
\begin{displaymath}
u_k(x,t)=\frac{u(d_kx+x_k,d_k^2t+t_k)}{d_k}, \chi_k(x,t)=\chi(d_kx+x_k,d_k^2t+t_k).
\end{displaymath}
We know that, passing to a subsequence if necessary,
$u_k\rightarrow u_0$ locally uniformly in $\R^{n+1}$
and $\chi_k\rightharpoonup \chi_0$ weakly-* in $L^\infty_{\rm loc}(\R^{n+1})$ as $k\to\infty$.
Also, after
a rotation and translation, the scaled
$B$ converges to $\{ x_n > 0 \}$ and $\left((y_k-x_k)/d_k, (s_k-t_k)/d_k^2\right)\rightarrow (\xi,\tau)
\in \partial Q_1(0)$ as $k\to \infty$.
The limit function $u_0$ satisfies
\begin{displaymath}
\begin{array}{ll}
\Delta u_0 -\partial_t u_0= 0 & \textrm{in } \R^{n+1}\cap 
\{ u_0>0 \} \; ,\\
u_0(\xi,\tau)=\ell & \textrm{ and}\\
u_0(x,t)=0 & \textrm{in } \{x_n > 0\}.
\end{array}
\end{displaymath}
By the definition of the limit superior we know also that
\begin{displaymath}
u_0(x,t)\le -\ell x_n \quad \textrm{ in } \{x_n<0\}.
\end{displaymath}
The strong maximum principle (applied to $u_0(x,t)+\ell x_n$) tells us therefore
that $u_0(x,t)=\ell \max(-x_n,0)$ for $t<\tau$. We have to
show that $\ell = 1$.
\\
In the case $\ell>0$ we obtain from 
the fact that $(u,\chi)$ is a solution in the sense of
Definition \ref{solution}, that $\chi_0 =1$ in 
$\{ x_n<0\}\cap \{ t<\tau\}$. Furthermore, we infer from the assumption that
$\chi_0 =0$ in $\{ x_n>0\}$. But then 
$(u_0(\cdot,t+\tau),\chi_0(\cdot,t+\tau))$ is in $\{ t<0\}$ a solution in the sense of
Definition \ref{solution} whose
energy
$M(u_0(\cdot,t+\tau),\chi_0(\cdot,t+\tau))=H_n$ (cf. \cite[Section 10]{calc}), whence \cite[Proposition 10.1]{calc}
implies that $\ell=1$.
\\
In the case $\ell=0$ we apply Lemma \ref{prem} to obtain
for some $\delta>0$ that $u$
is caloric in $ Q_\delta$ and satisfies
\[ u=0 \textrm{ in } Q^-_\delta\; .\]
As $\{ u= 0\}$ contains $B$, $u$ being caloric in $Q_\delta$
and therefore analytic with respect to the space variables implies
\[ u=0 \textrm{ in } Q_{\delta_1}\]
for some  $\delta_1>0$.
This is a contradiction in view of the origin being a free boundary point.
\qed
\\
The following theorem extends \cite[Lemma 7.2]{AC}.
\begin{theorem} \label{lem:flatabove}
There exists a constant $C\in(0,+\infty)$ 
depending only on the space dimension $n$
such that if 
$u\in F(\sigma,1,\sigma)$ in $Q_{\rho}$ then $u\in F(C\sigma,C\sigma,\sigma)$
in $Q_{\rho/2}(0,y_n,0))$ for some $|y_n|\le C\sigma$.
\end{theorem}
\proof
The idea is to touch the boundary $\partial \{ \chi=0\}$ with the graph of a $C^2$-function, 
to apply Lemma \ref{lem:ball} and to proceed then with a Harnack inequality
argument.\\
{\bf Step 1 (Touching $\partial \{ \chi=0\}$ with a smooth surface):}
\\
Rescaling 
$u_{\rho}(x,t) := {u({\rho}x,{\rho}^2 t)\over {\rho}}, \chi_{\rho}(x,t) := \chi({\rho}x,{\rho}^2 t)$
we may assume that $\rho=1$.
Let 
\begin{displaymath}
\eta(x',t)=\left\{\begin{array}{ll}
\exp(\frac{16(|x'|^2+|t-1|)}{1-16(|x'|^2+|t-1|)}),& |x'|^2+|t-1|<1/16,\\
0,&  \textrm{ else}\end{array}\right.
\end{displaymath}
and let $s$ be the largest constant such that
\begin{displaymath}
Q_1\cap \{ u>0\} \subset \{ (x,t)\in Q_1:\; x_n<\sigma-s\eta(x',t) \}=: D.
\end{displaymath}
This implies that there exists a point
$(x_0,t_0):=Z\in \partial D\cap \partial \{ u>0\}\cap \{ t\ge 15/16\}$.
As $(0,1)$ is a free boundary point,
we know furthermore that $s\le \sigma$.
\\
Let us also define the barrier function $v$ by
\begin{displaymath}
\begin{array}{ll}
\Delta v - \partial_t v =0 & \textrm{ in } D \; ,\\
v=0 & \textrm{ on } \partial D \cap Q_1 \textrm{ and}\\
v=2\sigma -x_n & \textrm{ on }  \partial D\cap\partial'Q_1.
\end{array}
\end{displaymath}
Note that this implies that $-\sigma\le v+x_n\le 2\sigma$.
\\
Since $|\nabla u|\le 1+\sigma$ we also obtain that $v\ge u$ on $\partial D$ 
and thus, by the maximum principle, that $v\ge u$ in $D$. As $\eta$
is a $C^2$-function, the assumptions of Lemma \ref{lem:ball} are satisfied at $Z$.
Therefore
\begin{equation}\label{eq:norest}
1\le
\limsup_{(x,t)\rightarrow Z}\frac{u(x,t)}{\pardist((x,t),B)}\le
-\partial_\nu v(Z),
\end{equation}
where $\nu$ is the outward space normal to $\partial D$ at $Z$. 
In order to obtain an estimate from above we define
\begin{displaymath}
F(x,t)=2\sigma-x_n-v(x,t).
\end{displaymath}
$F$ is caloric in $D$ and satisfies $0\le F\le \sigma$. Since $D$ is a regular parabolic
domain, we know from standard regularity theory for
parabolic equations that $\sup_D |\nabla F|\le C_1\sigma$. Therefore
\begin{displaymath}
-\partial_n v(Z)=1+\partial_n F(Z)\le
1+C_1\sigma.
\end{displaymath}
By the flatness assumption
we know that  $\nu$ is close to $e_n$. More precisely,
\begin{displaymath}
|\nu-e_n|=\big| 
\frac{(-s\nabla \eta,1-\sqrt{s^2|\nabla \eta|^2+1})}{\sqrt{s^2|\nabla \eta|^2+1}}\big|
\le \sqrt{10} |\nabla \eta| s.
\end{displaymath}
Thus
\begin{displaymath}
-\partial_\nu v(Z)=-\nabla v(Z)\cdot(\nu-e_n) -
\partial_n v(Z) \le 1+C_1\sigma +\sqrt{10} |\nabla \eta| |\nabla v(Z)| s \le 1+ C_2 \sigma\; .
\end{displaymath}
From inequality (\ref{eq:norest}) we infer that
\begin{equation}\label{eq:bar}
1\le -\partial_\nu v(Z) \le 1+C_2 \sigma.
\end{equation}
\\
{\bf Step 2 (Harnack inequality argument):}\\
As we know already that $v$ is $\sigma$-close to $-x_n$, it is
sufficient to show that 
$u$ is $\sigma$-close to $v$ on the set
$\{(x,t):\; x_n=-3/4,\; |x'|\le 1/2,\; t\le 3/4\}$. Once this is done,
we may integrate $u$
in the $x_n$-direction to establish the lemma.
\\
In order to prove the $\sigma$-closeness
we define 
for $\xi=(\gamma,\tau)$, $\tau\in (-1,3/4)$, $|\gamma'|\le 1/2$ and 
$\gamma_n=-3/4$
the function $\omega_\xi$ by
\begin{displaymath}
\begin{array}{ll}
\Delta \omega_\xi - \partial_t \omega_\xi= 0 &
\textrm{ in } D\cap \{t>\tau\} \\
\omega_\xi = -x_n & \textrm{ on } B_{1/8}(\gamma)\times \{t=\tau\} \\
\omega_\xi = 0 & \textrm{ on the remainder of the parabolic boundary of } 
D\cap \{t>\tau\}.
\end{array}
\end{displaymath}
By the Hopf lemma we have 
\begin{displaymath}
\partial_\nu \omega_\xi(Z)\le -\alpha <0
\end{displaymath}
uniformly in $\xi$.
\\
We would like to show that $u\ge v-C_4\sigma x_n$. The trick is to compare $u$
to $v-K\sigma \omega_\xi$ on the set $B_{1/8}(\gamma)\times \{t=\tau\}$ and to use 
the information on the normal derivative of $u$ at $Z$ to prove that if $K$
is large, then  $u>v-K\sigma \omega_\xi$ for at least one point in 
$B_{1/8}(\gamma)\times \{t=\tau\}$. More precisely:
\\
Assume that $u\le v-K\sigma \omega_\xi$ in 
$B_{1/8}(\gamma)\times \{t=\tau\}$. Then
$u\le v-K\sigma \omega_\xi$ in $D\cap \{t>\tau\}$. Consequently, we obtain from inequalities
(\ref{eq:norest}) and
(\ref{eq:bar}) that
\begin{displaymath}
1\le -\partial_\nu v(Z)+
K\sigma\partial_\nu \omega_\xi(Z)\le 1+
C_2\sigma-K\alpha\sigma.
\end{displaymath}
This yields a contradiction when $K$ is large enough, say $K = 2C_2/\alpha$. Thus
$u(X_\xi)>v(X_\xi)-K\sigma\omega_\xi(X_\xi)$ for at least one point  
$X_\xi\in B_{1/8}(\gamma)\times \{t=\tau\}$.
\\
On the other hand, $v-u\ge 0$. Therefore we can apply the Harnack inequality and deduce
that
\begin{displaymath}
(v-u)(\tilde\xi)\le C_3\inf_{Q_{1/8}(\tilde\xi+(0,1/32))}(v-u)\le C_4\sigma,
\end{displaymath}
for every $\tilde\xi\in \{(x',-3/4,t):\; |x'|<1/2,\; -1\le t\le 1/2 \}$.
\\
This implies that $u(x',-3/4,t)\ge 3/4 - C_5\sigma$ in the above region. 
Integrating in the $e_n$ direction and using the assumption $|\nabla u|\le 1+\sigma$ 
yields the estimate
\begin{displaymath}
u\ge -(x_n+C_6\sigma) \textrm{ in } \{ -3/4 \le x_n \le -\sigma\} \times Q'_{1/2}
\end{displaymath}
By our
initial assumption we also know that $u=0$ in $\{3/4 \ge x_n\ge \sigma\}\cap Q'_{1/2}$.
Translating $(u,\chi)$
in the $e_n$ direction so that the point $(0,1/4)\in \partial\{ u>0\}$ 
and using $\chi\ge \chi_{\{ u>0\}}$ of Definition \ref{solution} 1)
we obtain the statement of our
theorem.
\qed
\section{Inhomogeneous Blow-up}
In this section we consider inhomogeneous scaling of the solution
and the free boundary. The following lemma is our version
of \cite[Lemma 7.3]{AC}
\begin{lemma}\label{lem:nonhomogeneousblowup}
Suppose that  $u_k \in F(\sigma_k,\sigma_k,\tau_k)$ in $Q_{\rho_k}$, 
that $\sigma_k\rightarrow 0$ and that $\tau_k/\sigma_k^2\rightarrow 0$, and
define
\begin{displaymath}
f_k^+(x',t):=\sup\{h:\; \limsup_{r\to 0}r^{-n-2}\int_{Q_r(\rho_k x', \sigma_k\rho_k h, \rho_k^2 t)}\chi > 0\},
\end{displaymath}
\begin{displaymath}
f_k^-(x',t):=\inf\{h:\; \limsup_{r\to 0}r^{-n-2}\int_{Q_r(\rho_k x', \sigma_k\rho_k h, \rho_k^2 t)} \chi > 0\}.
\end{displaymath}
Then, as a subsequence $k\to\infty$,
$f_k^+$ and $f_k^-$ converge in $L^\infty_{\rm loc}(Q'_1)$ to some function $f$,
and $f$ is continuous in $Q'_1$.
\end{lemma}
\proof
Rescaling as before we may assume that $\rho_k=1$. Let 
\[ D_k:=\{ (y',h,t):\; \limsup_{r\to 0}r^{-n-2}\int_{Q_r(y',\sigma_k h, t)} \chi > 0 \}\; .\] 
We may assume -- passing if necessary to a subsequence -- 
that $D_k$ converges with respect to the usual (not the parabolic) Hausdorff distance
as $k\to\infty$.
Let us define \begin{displaymath}
f(x',t):=  
\limsup_{(y',s)\rightarrow (x',t), k\rightarrow \infty} f_k^+(y',s),
\end{displaymath}
where we take the limit superior with respect to the above subsequence.
For every $(y_0',t_0)$ there exists then a sequence
$(y'_k,t_k)\rightarrow (y'_0,t_0)$ such that 
$f^+_k(y_k',t_k)\rightarrow f(y'_0,t_0)$ as $k\to\infty$.
By definition $f$ is upper semi-continuous. Therefore we obtain for
 $\varepsilon>0$ and sufficiently large $k$ that
\begin{displaymath}
\big( \overline{Q_{\varepsilon}'(y_k',t_k)}\times [f_k^+(y'_k,t_k)+\delta,\infty)\big) \cap \bar{D}_k=\emptyset. 
\end{displaymath}
\\
Consequently $u_k\in F(\sigma_k\frac{\delta}{\varepsilon},1,\tau_k)$ 
in $Q_{\varepsilon}(y_k,\sigma_k f^+_k(y'_k,t_k),t_k)$. Applying Theorem
\ref{lem:flatabove} to $u_k$ we deduce that 
\[u_k(x,t) \ge -(x_n +C\sigma_k\delta/2) \textrm{ for } (x,t)\in Q_{\varepsilon/2}(y_k',\sigma_k f^+_k(y'_k,t_k),t_k)\; .\]
In terms of $f^+_k$ and $f^-_k$ this yields
$f^-_k(y',t)\ge f^+_k(y'_k,t_k)-C\delta$ in $Q'_{\varepsilon/4}(y_k',t_k)$.
It follows that
$\lim_{k\to\infty}f^-_k(y',t)=f(y',t)$, that
$f^+_k$ and $f^-_k$ converge locally uniformly and that
$f$ is continuous.
\qed
\\
The next Proposition follows the lines of \cite[Lemma 5.7]{ACF}.
\begin{proposition}\label{lem:wcal} Suppose that the assumptions
of Lemma \ref{lem:nonhomogeneousblowup} are satisfied and that $k$ is the
subsequence of Lemma \ref{lem:nonhomogeneousblowup}.
Then
\begin{displaymath}
w_k(x',h,t)=\frac{u_k(\rho_k x',\rho_k h,\rho_k^2 t)+\rho_k h}{\sigma_k}
\end{displaymath}
is for each $\delta\in (0,1)$ bounded in
$Q_{1-\delta}\cap\{ x_n< 0\}$ (by a constant depending only on $\delta$ and $n$)
and
converges on compact subsets of $Q_1^-$ in $C^2$ to a caloric function $w$.\\
Moreover, 
$w(x',h,t)$ is non-decreasing in the $h$-variable in
$Q_1^-$ and
\[ \lim_{Q_1^- \ni (y,s)\to (x',0,t)\in Q'_1,k\to\infty} w_k(y,s)=f(x',t)\> ;\]
here $f$ is the function defined
in Lemma \ref{lem:nonhomogeneousblowup}.
\end{proposition}
\proof
Rescaling as before we may assume that $\rho_k=1$.\\
The function $w_k$ is caloric in 
$Q_1\cap \{h<-\sigma_k\}$. 
Using Definition \ref{flatness} 3), we obtain that
\[ u_k\le -x_n +2\sigma_k \textrm{ in } Q_1\cap \{x_n \le 0\}\; ,\] 
implying that $w_k\le 2$. From Theorem \ref{lem:flatabove} 
and
Definition \ref{flatness} 3)
we infer
furthermore that
$u_k(x,t)\ge -(x_n + C_\delta\sigma_k)$ for 
$(x',x_n,t)\in Q_{1-\delta}\cap \{ x_n\le 0\}$,
implying that $w_k \ge -C_\delta$ in $Q_{1-\delta}\cap \{ x_n\le 0\}$.
\\
By Definition \ref{flatness} 3) and the assumptions,
$|\nabla u_k|\le 1+o(\sigma_k^2)$. Consequently,
\begin{equation}\label{eq:hder}
-\partial_h w_k\le\frac{|\nabla u_k|-1}{\sigma_k}\le
\frac{\tau_k}{\sigma_k}\rightarrow 0\textrm{ as } k\to\infty\; .
\end{equation}
In the remainder of the proof we will show that $w$ attains the boundary data
$f$ as $h\to 0$.
First, we show that for fixed $L\in (1,+\infty)$
\begin{equation}\label{eq:interi}
w_k(x',\sigma_k h,t)-f_k^+(x',t)\rightarrow 0 \quad 
\textrm{ uniformly in } Q'_{1-\delta}\times \{ -L\le h<0 \}
\end{equation}
as $k\to\infty$.
An estimate from above can be obtained easily from 
inequality (\ref{eq:hder}):
\begin{displaymath}
w_k(x',h\sigma_k,t)-f^+_k(x',t)\le w_k(x',\sigma_k f^+_k(x',t),t)-f^+_k(x',t)+
(f^+_k(x',t)-h)\frac{\tau_k}{\sigma_k}
\end{displaymath}
\begin{displaymath}
\le (1+L)\frac{\tau_k}{\sigma_k}
\rightarrow 0 \textrm{ as } k\to\infty \; .
\end{displaymath}
This establishes an estimate from above. In order to derive an estimate
from  below we use Theorem \ref{lem:flatabove}:
Consider a sequence of points $(x'_k,t_k)\in Q_{1-\delta}'$ and fixed $S\in (4,+\infty)$.
Then
\begin{displaymath}
u_k\in F(\tilde{\sigma}_k,1,\tau_k) 
\textrm{ in } Q_{S\sigma_k}(x'_k,\sigma_k f_k^+(x'_k,t_k),t_k)
\end{displaymath}
for
\begin{displaymath}
\tilde{\sigma}_k=
\frac{1}{S}\sup_{(x',t)\in Q_{S\sigma_k}'}(f^+_k(x',t)-f^+_k(x'_k,t_k)).
\end{displaymath}
From the uniform convergence of $f^+_k$ to the continuous function
$f$, we infer that 
$\tilde{\sigma}_k\rightarrow 0$ as $k\to\infty$. Now by Theorem \ref{lem:flatabove},
\begin{displaymath}
u_k\in F(C\bar{\sigma}_k,C\bar{\sigma}_k,\tau_k) \quad 
\textrm{ in } Q_{S\sigma_k/2}(x'_k,\sigma_k f_k^+(x'_k,t_k)+CS\bar\sigma_k\theta/2,t_k),
\end{displaymath}
where $\bar{\sigma}_k=\max(\tilde{\sigma}_k, \tau_k)$ and $\theta \in [0,1]$.
\\
Thus for $h\in (\max(-L,-S/4),0)$
\begin{displaymath}
u_k(x_k+h\sigma_ke_n,t_k)\ge 
-\sigma_k \left( h-f_k^+(x'_k,t_k)+C\bar\sigma_kS\right).
\end{displaymath}
Consequently
\begin{displaymath}
w_k(x_k+h\sigma_ke_n,t_k)=
\frac{u_k(x_k+h\sigma_ke_n,t_k)+h\sigma_k}{\sigma_k}\ge 
f_k^+(x'_k,t_k)-C\bar{\sigma}_kS\; ,
\end{displaymath}
and (\ref{eq:interi}) holds.
\\
To establish  $\lim_{Q_1^- \ni (y,s)\to (x',0,t)\in Q'_1,k\to\infty} w_k(y,s)=f(x',t)$, we need to extend the
convergence (\ref{eq:interi}) to larger values of $h$. To this end,
we define the barrier function $z_\varepsilon$ by
\begin{displaymath}
\begin{array}{ll}
\Delta z_\varepsilon -\partial_t z_\varepsilon = 0 &
\textrm{ in } Q^-_{1-\delta} \; ,\\
z_\varepsilon =g_\varepsilon & \textrm{ on } \partial' Q_{1-\delta}\cap \{ h=0 \}\; ,\\
z_\varepsilon = \inf_k \inf_{Q_{1-\delta}^-} w_k & \textrm{ on } 
\partial' Q_{1-\delta}\cap \{ h < 0 \}\; ,
\end{array}
\end{displaymath}
where  $g_\varepsilon \in C^{\infty}$ and 
$f-2\varepsilon\le g_\varepsilon\le f-\varepsilon$. By (\ref{eq:interi}) we know
that $w_k\ge z_\varepsilon$ on $\partial'(Q_{1-\delta}\cap \{h\le -L\sigma_k\})$.
From the comparison principle it follows that $w_k\ge z_\varepsilon$ in 
$Q_{1-\delta}^-\cap \{h\le -L\sigma_k\}$. Thus, by local boundary regularity for solutions of
the heat equation, 
$\liminf_{Q_{1-2\delta}^- \ni (y,s)\to (x',0,t),k\to\infty} w_k(y,s)\ge g_\varepsilon(x',t)\ge f(x',t)-2\varepsilon$.\\
The opposite inequality follows from a similar argument, comparing $w_k$
to the upper barrier $\tilde z$ defined by
\begin{displaymath}
\begin{array}{ll}
\Delta \tilde z_\varepsilon -\partial_t \tilde z_\varepsilon = 0 &
\textrm{ in } Q^-_{1-\delta} \; ,\\
\tilde z_\varepsilon =\tilde g_\varepsilon & \textrm{ on } \partial' Q_{1-\delta}\cap \{ h=0 \}\; ,\\
\tilde z_\varepsilon = \sup_k \sup_{Q_{1-\delta}^-} w_k & \textrm{ on } 
\partial' Q_{1-\delta}\cap \{ h < 0 \}\; ,
\end{array}
\end{displaymath}
where  $\tilde g_\varepsilon \in C^{\infty}$ and 
$f+2\varepsilon\ge \tilde g_\varepsilon\ge f+\varepsilon$.
\qed
\section{Scaling discrepancy and $C^\infty$-regularity of blow-up limits}\label{discrep}
In order to obtain ``better-than-Lipschitz''-regularity of the inhomogeneous blow-up limit $f$, 
H.W. Alt-L.A. Caffarelli used the non-positive
mean curvature of $\partial\{ u>0\}$ at singularities. 
The analogue of the non-positive
mean curvature property can still be proved in the time-dependent case, 
however that path leads to problems in the sequel. Therefore
we replace it by a scaling discrepancy argument which gives hope to be applicable
in more general situations. We obtain $C^\infty$-regularity of $f$.
\begin{proposition}\label{discr}
Suppose that the assumptions
of Lemma \ref{lem:nonhomogeneousblowup} are satisfied and that $k$ is the
subsequence of Lemma \ref{lem:nonhomogeneousblowup}. Then
$\partial_n w=0$ on $Q'_{1/2}$ in the sense of distributions.
\end{proposition}
\proof
Rescaling as before we may assume that $\rho_k=1$.\\
In what follows, $g(x',t)=8(|x'|^2+ |t|)-4$. Note that $f\ge g$ in $Q'_{1/2}$.
Let us introduce the following notation:
$Z$ shall be the set $\{ (x',x_n,t): (x',t)\in Q'_1,x_n\in \R\}$. 
Given a function $\phi: Q'_1\to \R$, we divide $Z$ into the three parts
\begin{displaymath}
Z^+(\phi)=\{ (x,t)\in Z: \; x_n> \phi(x',t) \},
\end{displaymath}
\begin{displaymath}
Z^-(\phi)=\{ (x,t)\in Z: \; x_n< \phi(x',t) \},
\end{displaymath}
\begin{displaymath}
Z^0(\phi)=\{ (x,t)\in Z: \; x_n = \phi(x',t) \}.
\end{displaymath}
Moreover let $\mu$ be defined by $\mu(A) := \int_{-\infty}^\infty {\mathcal H}^{n-1}(A \cap \{ s=t \})\> dt$
for any Borel set $A\subset \R^{n+1}$.
Adding an arbitrarily small constant to the function $g$, we may assume
that $\mu(Z^0(\sigma_k g)\cap R_k)=0$
for all $k$; here $R_k$ is the regular part of the
free boundary $\partial\{ u_k>0\}$ introduced in \cite[Proposition 9.1]{calc}, i.e.
$$R_k(t) :=\{ x\in \partial\{u_k(t)>0\}\> : \>
\textrm{there is } \nu_{R_k}(x,t)\in \partial B_1(0) \textrm{ such that }
v_r(y,s) = $$ $${u_k(x+ry,t+r^2s)\over r} \> \to \>
\max(-y\cdot \nu_{R_k}(x,t),0) \textrm{ locally uniformly in } (y,s)\in
\R^{n+1} $$ $$\textrm{ as } r \to 0\}\; .$$
Last, we define $E_k := \{ u_k >0\}\cap Z^-(\sigma_k g)$
and $\Sigma_k := \{ (x',t) : (x',\sigma_k g(x',t),t)\in \{ u_k>0\}\cap Z\}$.
By the choice of $g$ we know that the limit inferior of the sets $\Sigma_k$
contains $Q'_{1/2}$.
\\
We will deduce the result from the following three claims.\\
\textsl{Claim 1:}
\begin{displaymath}
\mu(Z^+(\sigma_k g)\cap  {R_k})\le -\int_{{\Sigma_k}}  (\partial_n u_k+1)dx'dt+
{\mathcal L}^n({\Sigma_k})+C_1\sigma_k^2.
\end{displaymath}
\textsl{Claim 2:}
$$
{\mathcal L}^n({\Sigma_k})-C_2\sigma_k^2 \le
\mu(Z^+(\sigma_k g)\cap R_k).
$$
\textsl{Claim 3:} 
$$
\int_{{\Sigma_k}}|\partial_n w_k(x',\sigma_k g(x',t),t)| \; \to \; 0 \textrm{ as } k\to\infty\; .
$$
\textsl{Proof of Claim 1:} By the representation theorem
\cite[Lemma 11.3]{calc}
we know that for non-negative $\phi\in C^{\infty}_0$,
\begin{equation}\label{eq:claim1}
\int_{-\infty}^{\infty}\int_{{R_k}(t)}\phi \> d{\mathcal H}^{n-1}\> dt \le
-\int_{\{ u_k>0\}}\left(\nabla u_k \cdot \nabla \phi +
\partial_t u_k\phi \right)\> dx\> dt\; .
\end{equation}
Letting $\phi\rightarrow \chi_{Z^+(\sigma_k g)}\chi_{Q_2}$ the inequality (\ref{eq:claim1})
becomes
\begin{equation}\label{eq:troll}
\mu(Z^+(\sigma_k g)\cap  {R_k})=
\int_{-\infty}^{\infty}\int_{{R_k}(t)\cap Z^+(\sigma_k g)} d{\mathcal H}^{n-1}\> dt
\end{equation}
\begin{displaymath}
\le \;
\int_{\{ u_k>0\}\cap Z^0(\sigma_k g)}\nabla u_k \cdot \nu \> dx\> dt -
\int_{\{ u_k>0\}\cap Z^+(\sigma_k g)}\partial_t u_k \> dx\> dt\; ,
\end{displaymath}
where $\nu$ is the outward unit space normal on $\partial Z^+(\sigma_k g)$. 
Since 
$$
\nu=\frac{1}{\sqrt{1+|\sigma_k \nabla' g|^2}}(\sigma_k\nabla' g,-1)\; ,
$$
we obtain
$$\mu(Z^+(\sigma_k g)\cap  {R_k})\;
\le\;
\int_{{\Sigma_k}}(\nabla u_k)(x',\sigma_k g(x',t),t) \cdot (\sigma_k \nabla' g(x',t), -1) \> dx\> dt$$ $$-\> 
\int_{\{ u_k>0\}\cap Z^+(\sigma_k g)}\partial_t u_k \> dx\> dt\; .
$$
Let us rewrite the integral
$$
\int_{{\Sigma_k}}(\nabla u_k)(x',\sigma_k g(x',t),t) \cdot (\sigma_k \nabla' g(x',t), -1) \> dx\> dt
$$
$$ = \;
\int_{{\Sigma_k}}\sigma_k(\nabla' u_k)(x',\sigma_k g(x',t),t) \cdot \nabla'g(x',t) - 
(\partial_n u_k(x',\sigma_k g(x',t),t)+1)\> dx'\> dt \; + \; {\mathcal L}^n({\Sigma_k})
$$
$$
= \int_{{\Sigma_k}} -\sigma_k u_k(x',\sigma_k g(x',t),t) \Delta' g(x',t)
\> - \>{\sigma_k}^2 \partial_n u_k(x',\sigma_k g(x',t),t) |\nabla' g(x',t)|^2
$$
$$
- \> (\partial_n u_k(x',\sigma_k g(x',t),t)+1)\> dx'\> dt + {\mathcal L}^n({\Sigma_k})
$$
$$
+\> \int_{\partial {\Sigma_k}} \sigma_k u_k(x',\sigma_k g(x',t),t) \partial_\eta g(x',t)\> d{\mathcal H}^{n-2} \> dt\; ,
$$
where $\eta$ is the outward space normal on $\partial {\Sigma_k}$.
Since $u_k=0$ on $\partial {\Sigma_k}$, the last integral is $0$.
\\
Moreover, $\Delta' g= 16$ and $u_k \le C_3 \sigma_k$ on $(x',g(x',t),t)$,
implying that
$$
\int_{{\Sigma_k}}(\nabla u_k)(x',\sigma_k g(x',t),t) \cdot (\sigma_k \nabla' g(x',t), -1) \> dx\> dt
$$
$$ =\>
-\int_{{\Sigma_k}}  (\partial_n u_k(x',\sigma_k g(x',t),t)+1)dx'dt+
{\mathcal L}^n({\Sigma_k})+C_4\sigma_k^2\; .
$$
By the definition of $w_k$ this tells us also that
\begin{equation}\label{later}
\int_{{\Sigma_k}}(\nabla w_k)(x',\sigma_k g(x',t),t) \cdot (\sigma_k \nabla' g(x',t), 0) \> dx\> dt
\; \to\; 0 \textrm{ as } k\to\infty\; ,
\end{equation}
a fact that will be used later on.\\
Last, integration by parts of the last term in (\ref{eq:troll})
with respect to the time variable yields
$$-
\int_{\{ u_k>0\}\cap Z^+(\sigma_k g)}\partial_t u_k \> dx\> dt\;
 \le \; 
C_5 {\sigma_k}^2\; .
$$ 
Combining the above estimates we obtain Claim 1.\\
\textsl{Proof of claim 2:}
With the outward space normal on the boundary of $Z^-(\sigma_k g)$
\begin{displaymath}
\nu_{g_k}=
\frac{1}{\sqrt{1+\sigma_k^2|\nabla' g|^2}}(-\sigma_k\nabla' g,1)
\end{displaymath}
and 
with the outward space normal $\nu_{R_k}$ on the regular boundary of $E_k$ we compute
\begin{equation}\label{eq:eta}
\mu(Z^+(\sigma_k g) \cap R_k)\ge 
\end{equation}
\begin{displaymath}
\int_{-1}^{1}\int_{Z^+(\sigma_k g)\cap 
R_k(t)}
\nu_{g_k}\cdot \nu_{R_k}\> d{\mathcal H}^{n-1}\> dt
\end{displaymath}
\begin{displaymath}
=\; \int_{-1}^1 \int_{E_k\cap Z^+(\sigma_k g)} \div\nu_{g_k}\> d{\mathcal H}^{n-1}\> dt
\end{displaymath}
\begin{displaymath}
+\; 
\int_{-1}^1\int_{\partial Z^+(\sigma_k g)\cap E_k} \nu_{g_k}\cdot 
\nu_{g_k}d{\mathcal H}^{n-1}\> dt.
\end{displaymath}
The normal $\nu_{g_{k}}$ satisfies
\begin{displaymath}
\div \nu_{g_k}\ge
\frac{-\sigma_k \Delta g}{\sqrt{1+\sigma_k^2|\nabla' g|^2}}
\ge - C_6\sigma_k. 
\end{displaymath}
Inserting this estimate for the divergence into (\ref{eq:eta}) yields
\begin{equation}\label{eq:est1}
\mu(Z^+(\sigma_k g)\cap R_k)\ge
\mu(\partial Z^+(\sigma_k g)\cap E_k)
\end{equation}
\begin{displaymath}-\; 
\int_{-1}^1 \int_{E_k\cap Z^+(\sigma_k g)} C_6 \sigma_k  \> d{\mathcal H}^{n-1}\> dt
\end{displaymath}
\begin{displaymath}
\ge \; \mu(\partial Z^+(\sigma_k g)\cap E_k)\> -\> C_7 \sigma_k^2\; ;
\end{displaymath}
the last inequality follows from the 
fact that the width of the set $E_k$ is of order $O(\sigma_k)$.
As the area of $\partial Z^+(\sigma_k g)\cap E_k)$
is greater than that of $\Sigma_k$, the statement
of Claim 2 holds.\\
\textsl{Proof of Claim 3:} From \textsl{Claim 1} and \textsl{Claim 2} we 
infer that
$$
-C_8 \sigma_k^2 \le -\int_{{\Sigma_k}}  (\partial_n u_k(x',\sigma_k g(x',t),t)+1)\> dx'\> dt\; .
$$
But since $u_k\in F(\sigma_k,\sigma_k,\tau_k)$ and $\tau_k/\sigma_k^2\to 0$ 
as $k\to\infty$, it follows that
$$
\partial_n u_k+1\ge -|\nabla u_k|+1\;\ge \; -o(\sigma_k^2).
$$
Consequently
$$
\int_{{\Sigma_k}}|\partial_n w_k(x',\sigma_k g(x',t),t)|
\; = \; \int_{{\Sigma_k}}\left|\frac{\partial_n u_k(x',\sigma_k g(x',t),t)+1}{\sigma_k}\right|
$$
$$
\le \;
\int_{{\Sigma_k}}2\max\left(-\frac{\partial_n u_k(x',\sigma_k g(x',t),t)+1}{\sigma_k},0\right)
\> + \>
\int_{{\Sigma_k}}\frac{\partial_n u_k(x',\sigma_k g(x',t),t)+1}{\sigma_k}
$$
$$
\; \le
C_9 \sigma_k \to 0 \textrm{ as } k\to \infty\; ,
$$
and \textsl{Claim 3} is proved.\\
\textsl{Proof of the Proposition:}
Let $\zeta \in C^1_0(Q_{1/2})$. From Claim 3, from the fact that
$w_k$ is caloric in $Z^-(\sigma_k g)$, from (\ref{later})
and from a standard energy estimate for caloric functions
we infer now that
$$ o(1) \; = \; \int_{{\Sigma_k}}\zeta \>\partial_n w_k(x',\sigma_k g(x',t),t) \nu_n
$$
$$ = \; \int_{Z^-(\sigma_k g)} (\partial_n\zeta \> \partial_n w_k \> - \> \zeta \> \Delta' w_k\> + \> \zeta \> \partial_t w_k)
$$
$$
\; = \; o(1) \> + \> \int_{Z^-(\sigma_k g)} (\partial_n\zeta \> \partial_n w_k \> + \> \nabla'\zeta \cdot\nabla' w_k\> - \> w_k \> \partial_t \zeta)
$$
$$
\to \; \int_{Q^+_1} (\partial_n\zeta \> \partial_n w\> + \> \nabla'\zeta \cdot\nabla' w\> - \> w \> \partial_t \zeta)
\textrm{ as } k\to\infty\; ;
$$
here $\nu$ is the outward unit space normal on $\partial Z^-(\sigma_k g)$.
It follows that $\partial_n w=0$ on $Q'_{1/2}$ in the sense of distributions.
\qed
\begin{corollary}\label{cor:cinf} Suppose that the assumptions
of Lemma \ref{lem:nonhomogeneousblowup} are satisfied and that $k$ is the
subsequence of Lemma \ref{lem:nonhomogeneousblowup}. Then
$f\in C^{\infty}(Q_{1/2})$; 
moreover,
$$
\Big| \frac{\partial^{\alpha+k} f}{\partial x^\alpha\partial t^k}\Big|\le C(n,|\alpha|,k)
$$
in $Q_{1/4}$
for any $k\in \N$ and multi-index $\alpha \in \N^n$.
\end{corollary}
\proof
Since 
$\partial_n w=0$ on $Q'_{1/2}$ in the sense of distributions
we may reflect 
$w$ to a caloric function in $Q_{1/2}$.  
As $f=w|_{Q'_1}$
and $\Vert w \Vert_{L^\infty(Q_{3/4})}\le C(n)$
(see Proposition \ref{lem:wcal}),
the result follows from standard regularity theory
of caloric functions.\qed
\section{Flatness improvement and regularity}
Concluding regularity is then a standard procedure. 
The following Lemma \ref{lem:flatbelow}, Lemma \ref{lem:improve}
and Theorem \ref{theo:main} extend Lemma 7.9, Lemma 7.10 and
Theorem 8.1 in \cite{AC}.
Finally, we apply Theorem \ref{theo:main} to regular free boundary points,
i.e. points in the set $R$ 
defined in \cite[Proposition 9.1]{calc} (or the proof of Proposition \ref{discr})
to obtain that $R$ is open relative to $\partial \{ u>0\}$. 
\begin{lemma}\label{lem:flatbelow}
Let $\theta\in (0,1)$. Then
there exists a constant
$\sigma_{\theta}>0$ 
depending only on $\theta$ and the dimension $n$
such that if 
$\sigma<\sigma_{\theta}$, $\tau\le \sigma_{\theta}\sigma^2$ and
$u\in F(\sigma,\sigma,\tau)$
in $Q_{\rho}$ in direction $\eta$, then
\begin{displaymath}
u\in F(\theta\sigma,1,\tau) \; \textrm{ in } Q_{c(n) \theta \rho}(\vartheta\eta,0)
\end{displaymath}
in direction $\overline{\eta}$ for some $\vartheta\in [-\sigma,\sigma]$ and some $\overline{\eta}$
satisfying
$|\overline{\eta}-\eta|\le {C(n)\sigma}$.
Here $c(n)>0$ and $C(n)<+\infty$ are constants depending only
on the dimension $n$. 
\end{lemma}
\proof
We may rotate the coordinate system so $\eta=e_n$, and we may
assume that $\rho=1$.
By a contradiction argument, it is sufficient to prove the statement
of the lemma for $u_k$ as in Lemma 
\ref{lem:nonhomogeneousblowup} and every large $k$.
\\
First, observe that by Corollary \ref{cor:cinf},
$$
f(x',t)\le f(0,0)+\ell\cdot x' + C(|x'|^2+|t|) \textrm{ in } Q'_{1/4},
$$
where $\ell$ is the space gradient of $f$, $|\ell|\le C$
and $C$ depends only on the dimension $n$.
Thus
$$
f(x',t)\le f(0,0)+\ell\cdot x' +\frac{\theta }{4}\frac{\theta}{4C} \textrm{ in }
Q_{\theta/(4C)}\; .
$$
It follows that for large $k$ the function $f^+_k$ in 
Lemma \ref{lem:nonhomogeneousblowup} satisfies 
\begin{displaymath}
f^+_k(x',t)\le f(0,0) + \ell\cdot x' +\theta \frac{\theta}{4C} \quad \textrm{ in } Q_{\theta/(4C)}\; .
\end{displaymath}
This means that $u_k \in F(\theta \sigma,1,\tau)$ in $Q_{\theta/(4C)}(0,f(0,0),0)$ in the 
direction 
$\bar{\eta}$, where
\begin{displaymath}
\bar{\eta}=\frac{(-\sigma_k\ell,1)}{\sqrt{1+|\sigma_k \ell|^2}}.
\end{displaymath}
The lemma follows. \qed
\begin{lemma}\label{gradient}
Let $u$ be a solution in the sense of Definition \ref{solution}. Then
$$ \max(|\nabla u|^2-1,0)(x,t) \to 0 \textrm{ as } 0< \pardist((x,t),\{ u=0\})\to 0\; .$$
\end{lemma}
\proof
Consider a sequence $\{ u>0\} \ni (x_k,t_k) \to (x_0,t_0)$ such that
$$1<\ell:=\limsup_{\{ u>0\} \ni (x,t)\to (x_0,t_0)} |\nabla u(x,t)|^2 = \lim_{k\to\infty} |\nabla u(x_k,t_k)|^2\; .$$
Setting $r_k := \pardist((x_k,t_k),\{ u=0\})$, the blow-up sequence
$$u_k(y,s) := {u(x_k+r_ky,t_k+{r_k}^2s) \over {r_k}},
\chi_k(y,s) := \chi(x_k+r_ky,t_k+{r_k}^2s)$$
converges to a solution
$(u_0,\chi_0)$ 
in the sense of Definition \ref{solution}
satisfying $u_0>0$ in $Q_1$, $|\nabla u_0|^2\le \ell$ 
and $|\nabla u_0(0)|^2=\ell$.
The strong maximum principle implies that
$u_0(y,s)=\ell \max(y\cdot e,0)$ in $\{ y\cdot e >0\}\cap \{ s<0\}$ for some $e\in \partial B_1$.
\\
From \cite[Theorem 11.1]{calc} we infer that
$$ \{ y\cdot e = 0\} \cap \{ s<0\}\subset \Sigma_{**}$$
up to a set of vanishing ${\mathcal L}^{n-1}$-measure,
where 
$$\Sigma_{**}(t) :=\{ x\in \partial\{u_0(t)>0\}\> : \>
\textrm{there is } \theta(x,t)\in (0,1] \textrm{ and }
\xi(x,t)\in \partial B_1(0) \textrm{ such}$$
$$\textrm{that }{u_0(x+ry,t+r^2s)\over r} \> \to \> \theta(x,t)
\vert y\cdot \xi(x,t)\vert \textrm{ locally uniformly}$$
$$\textrm{in } (y,s)\in
\R^{n+1} \textrm{ as } r \to 0\}\; .$$
However $\theta(x,t)\in (0,1]$ contradicts $\ell>1$.
\qed
\begin{lemma}\label{lem:improve} For every $\theta\in (0,1)$ there exist $\sigma_{\theta}>0$
and $c_\theta\in (0,1/2)$ depending only on $\theta$ and the dimension $n$
such that if $u\in F(\sigma,1,\tau)$ in $Q_{\rho}$ in direction
$\eta$ with $\sigma\le \sigma_{\theta}$ and $\tau\le \sigma_{\theta}\sigma^2$
then $u\in F(\theta \sigma,\theta\sigma,\theta^2\tau)$ in 
$Q_{c_\theta \rho}(\overline{y},0)$ in the direction $\overline{\eta}$ for some 
$\overline{y},\overline{\eta}$ satisfying
$|\overline{\eta}-\eta|\le C(n) \sigma$ and $|\overline{y}|\le C(n) \sigma$.
Here $C(n)$ depends only on the dimension $n$.
\end{lemma}
\proof
We may
assume that $\rho=1$.
\\
From Lemma \ref{lem:flatabove}  
we infer that $u\in F(C\sigma,C\sigma,\tau)$ in $Q_{1/2}(y,0)$ in direction
$\eta$
for some $y\in B_{C\sigma}$.
Consequently we may
apply Lemma \ref{lem:flatbelow} to deduce that 
for some $\theta_1$ to be determined later,
$u\in F(C\theta_1\sigma,1,\tau)$
in $Q_{c(n)\theta_1}(\tilde y,0)$ in the direction $\bar{\eta}$ such that
$|\eta-\bar{\eta}|\le C\sigma$ and $|\tilde y-\bar{y}|\le (C+1)\sigma<1/2$,
provided that $\sigma_{\theta}$ has been chosen small enough in terms
of $\theta_1$.
\\
In order to be able to continue we need to show improvement with respect to the $\tau$-variable.
To that end, observe that
$U=\max(|\nabla u|-1,0)$ is by Lemma \ref{gradient} a continuous subcaloric function in $Q_1$
with boundary values less than 
$\tau \chi_{\{ u>0\}}\le \tau \chi_{\{ x_n\le \sigma\}}$. We may therefore compare
$U$ to the caloric function with boundary values 
$\tau \chi_{\{ x_n\le \sigma\}}$. It follows that $0\le U \le (1-c_1)\tau$
in $Q_{1/2}$ for  some $c_1>0$ depending only on the dimension $n$.
Thus $u\in F(C\theta_1\sigma,1,(1-c_1)\tau)$ in $Q_{c(n)\theta_1}(\tilde y,0)$
in the direction $\bar{\eta}$. Choosing $\theta_0:=\sqrt{1-c_1}$
and $\theta_1 := \theta_0/C$ we obtain
$u\in F(\theta_0\sigma,1,\theta_0^2\tau)$ in $Q_{c_2\theta_0}(y,0)$ in the direction $\bar{\eta}$ such that
$|\eta-\bar{\eta}|\le C\sigma$, where $c_2\in (0,1)$ depends only on the dimension
$n$.
\\
Iterating this process we see that
$$
u\in F(\theta_0^m\sigma,1,\theta_0^{2m}\tau) \textrm{ in }
Q_{(c_2\theta_0)^m}(y_m,0)
\textrm{ in the direction }\bar{\eta}_m
$$
where  
$|\eta-\bar{\eta}_m|\le C(n) \sigma\sum_{j=0}^{m-1} \theta_0^j$
and $|y_m| \le C(n) \sigma\sum_{j=0}^{m-1} (c_2\theta_0)^j$.
\\
Applying once more Lemma \ref{lem:flatabove} 
and choosing $\theta_0 := \theta^{1\over m}/C$ we obtain the statement of the lemma.
\qed
\begin{theorem}\label{theo:main}
There exists a constant $\sigma_0>0$  such that
if $u\in F(\sigma,1,\tau)$ in $Q_\rho(t_0,x_0)$, $\sigma\le\sigma_0$ and $\tau \le \sigma_0\sigma^2$,
then the topological free boundary $\partial\{ u>0\}$ is in
$Q_{\rho/4}(t_0,x_0)$ the graph of a $\C^{1+\alpha,\alpha}$-function;
in particular the
space normal is H\"older continuous in $Q_{\rho/4}(t_0,x_0)$.
\end{theorem}
\proof
Using Lemma \ref{lem:improve} inductively we see that
\begin{equation}\label{eq:osccontroll}
\begin{array}{l}
u\in F(\theta^k\sigma,\theta^k\sigma,\theta^{2k}\tau) \textrm{ in }Q_{{c_{\theta \over 2}}^k\rho}(y,s)\textrm{ in the direction } \overline{\eta}^k \\
\textrm{where } |\overline{\eta}^k-\eta|\le C(n)\sigma \sum_{j=0}^{k-1}(2\theta)^j
\textrm{ and } |\overline{y}^k-y|\le C(n)\sigma\sum_{j=0}^{k-1}(2c_{\theta/2}\theta)^j\; ,
\end{array}
\end{equation}
provided that $(y,s)\in Q_{1/2}(t_0,x_0)\cap \partial\{ u>0\}$, $\theta < 1/4$ and $$\sigma_0 < \min(1/(4C(n)),\sigma_{\theta/2}/2)\; ;$$
here we sacrificed some flatness in order to keep the original free boundary
point $(y,s)$.
We obtain existence of the outward space normal $\nu$ on $Q_{1/2}(t_0,x_0)$.
Moreover, $\nu$ satisfies by
(\ref{eq:osccontroll})
$$
\textrm{osc}_{Q_{{c_{\theta/2}}^k\rho}(y,s)}\nu \le 
C(n,\theta) \theta^k \sigma\; ,
$$
which implies H\"older-continuity of $\nu$.\qed
\begin{corollary}\label{regular}
For each point $(x_0,t_0)$ of the set $R$, the topological free boundary $\partial\{ u>0\}$
is in an open neighborhood
of $(x_0,t_0)$ the graph of a $\C^{1+\alpha,\alpha}$-function;
in particular, the
space normal is H\"older continuous in an open space-time neighborhood of $(x_0,t_0)$.
\end{corollary}
\proof
The Corollary follows from \cite[Proposition 9.1]{calc} and Lemma \ref{gradient}.
 \qed
 \bibliographystyle{plain}
\bibliography{parber.bib}
\end{document}